\documentclass{birkmult} \usepackage{jms-bm}
\begin{document}

\def\secn#1{Section~\ref{sec:#1}}
\def\thm#1{Theorem~\ref{thm:#1}}
\def\lem#1{Lemma~\ref{lem:#1}}
\def\cor#1{Corollary~\ref{cor:#1}}
\def\prop#1{Proposition~\ref{pr:#1}}
\def\figr#1{Figure~\ref{fig:#1}}



\def\FTC/{{\small{FTC}}}

\newcommand{\abs}[1]{\lvert#1\rvert}
\newcommand{\R}{\mathbb{R}}
\newcommand{\Z}{\mathbb{Z}}
\renewcommand{\arc}[1]{\gamma_{#1}}
\newcommand{\len}[1]{\ell_{#1}}
\newcommand{\bdry}{\partial}
\newcommand{\after}{\circ}
\newcommand{\isom}{\cong}
\newcommand{\setm}{\smallsetminus}
\newcommand{\eps}{\varepsilon}
\newcommand{\lk}{\textrm{lk}}
\newcommand{\intr}{\textrm{int}} 
\newcommand{\half}{\tfrac12}
\newcommand{\arcsec}{\textrm{arcsec}}
\newcommand{\cutoff}{\pi/8}
\renewcommand{\d}{\partial}
\renewcommand{\phi}{\varphi}




\makeatletter
\newbox\overbox
\def\fakeover#1{\setbox\overbox\hbox{$#1$}\hbox
                         {$\overline{#1\hskip-\wd\overbox}$\hskip\wd\overbox}}
\def\overnoarrow#1{\mskip1.5mu\overline{\mskip-1.5mu#1}}
\def\overrightarrow#1{\mskip2mu\vbox{\m@th\ialign{##\crcr
   \rightarrowfill\crcr
   \noalign{\kern-.4pt               
        \kern-\fontdimen22\textfont2 
        \nointerlineskip}
   ${\mskip0mu\hfil\fakeover{#1}\hfil\mskip6mu}$\crcr}}\mskip-2mu}
\def\overleftrightarrow#1{\mskip-3mu\vbox{\m@th\ialign{##\crcr
   \leftarrowfill\hskip-.6em\rightarrowfill\crcr
   \noalign{\kern-.4pt               
        \kern-\fontdimen22\textfont2 
        \nointerlineskip}
   ${\mskip3mu\hfil\fakeover{#1}\hfil\mskip6mu}$\crcr}}\mskip-2mu}

\def\ray#1{{\smash{\overrightarrow{#1}}}}
\def\sline#1{{\smash{\overleftrightarrow{#1}}}}
\def\segment#1{{\smash{\overnoarrow{#1}}}}
\makeatother

\title[Convergence and Isotopy for FTC Graphs]{Convergence and Isotopy Type \\ for Graphs of Finite Total Curvature}
\date{March 12, 2003; revised \today}

\author{Elizabeth Denne}
\address{Department of Mathematics \& Statistics\\
Smith College\\ Northampton, MA 01063\\ USA}
\email{edenne@email.smith.edu}

\author{John M. Sullivan}
\address{\tuaddr32}
\email{sullivan@math.tu-berlin.de}



\begin{abstract}
Generalizing Milnor's result that an \FTC/ (finite total curvature) knot has
an isotopic inscribed polygon, we show that any two nearby knotted \FTC/ graphs
are isotopic by a small isotopy.  We also show how to obtain sharper
constants when the starting curve is smooth.  We apply our main
theorem to prove a limiting result for essential subarcs of a knot.
\end{abstract}

\maketitle

\section{Introduction}\index{isotopy|(}
A \emph{tame} knot or link type is one that can be represented
by a polygonal space curve.  It is clear that the corners
of a polygon can be rounded off, so that a \ix{tame link} can also be
represented by a smooth, even~$C^\infty$, curve.
Conversely, given a $C^2$ curve, it is relatively easy
to find an isotopic inscribed polygon.  With less smoothness,
the arguments become more intricate but are by now standard.
Crowell and Fox \cite[App.~I]{CFox} gave a detailed proof for $C^1$ curves,
while Milnor~\cite{milnor} introduced the class of
\emph{finite total curvature (\FTC/) curves},%
\index{finite total curvature (FTC)!curve}
and showed they also have isotopic inscribed polygons.

Thus any embedded \FTC/ curve represents a tame link type; but
tame links (even the unknot) can of course also be represented by
curves of infinite total curvature or even infinite length.
However, our work in \ix{geometric knot theory} has suggested
to us that \FTC/ links are a very useful class, which might be
considered as ``geometrically tame''.
(See~\cite{Sul-FTC}, in this volume, for a survey on \FTC/ curves.
We show for instance that this is the natural class for standard results
like Schur's comparison theorem.)

Here, with a view towards convergence results in geometric
knot theory, we examine the question of how close two space
curves must be in order to be (ambient) isotopic.  Of course,
the phenomenon of local knotting means that
our notion of distance between space curves
must control tangent directions as well as position.  (Otherwise, within
any distance~$\eps$ of any knot there would be infinitely many
composite knots---with extra summands contained in $\eps$-balls---and
even certain prime satellites of the original knot.)
We also show that, when the initial curves are close enough,
the ambient isotopy can be made arbitrarily small.

The arguments of~\cite{CFox,milnor} not only produce one
isotopic inscribed polygon, but in fact show that
all sufficiently fine inscribed polygons are isotopic to a
given $C^1$ or \FTC/ curve.  We recover this result,
since finely inscribed polygons are nearby in our sense
(with one caveat, discussed later).

We prove three versions of our theorem.  The first is restricted
curves of positive thickness, that is, to $C^{1,1}$ links, but
it allows us to get nearly optimal bounds (in terms of the thickness)
for how close the two curves must be before they are guaranteed
to be isotopic.
The second version, our main theorem, applies not just
to links but to arbitrary knotted \FTC/ graphs; it gives
us control over the distance points are moved by the isotopy.
The final version applies to $C^1$ links and
allows us to prove that the isotopy is small
even in a stronger $C^1$ sense.

As an application of the main theorem, we consider the notion
of essential arcs and secants of a knot, defined~\cite{dsDDS} in terms of
certain knotted $\Theta$-graphs.   We show that essential secants
remain essential in limits, which is important
for results~\cite{dsDenne} on quadrisecants.

\section{Definitions}

\begin{defn}
A \emph{(knotted) rectifiable graph}\index{knotted graph}
is a (multi-)graph---of fixed finite combinatorial type---%
embedded in space such that each edge
is a rectifiable arc.\index{rectifiable curve}  An \emph{\FTC/ graph}%
\index{finite total curvature (FTC)!graph}
is a rectifiable graph where each arc has
finite total curvature.
\end{defn}

We note that Taniyama~\cite{dsTaniyama-tcgraph} has considered
minimization of the total curvature of knotted graphs, including
analogs of the F\'ary/Milnor theorem.

For an introduction to the theory of knotted graphs, see~\cite{kauf89}:
in particular, we are interested in graphs with what Kauffman calls
topological (rather than rigid) vertices.  The $k$ edges incident
to a vertex~$v$ can be braided arbitrarily near~$v$ without
affecting the knot type of the graph.  We do not require
our graphs to be connected, so our results apply to
links as well as knots.

As mentioned above, in order to conclude that nearby
knots are isotopic, we need to use a notion of distance
that controls tangent vectors.  A standard notion would be,
for instance, $C^1$--convergence of $C^1$~curves.
Our goal, however, is to consider curves which need not be~$C^1$.
Recalling that a rectifiable curve has a well-defined
tangent vector almost everywhere, we can make the
following definition.

\begin{defn}
Given two rectifiable embeddings $\Gamma$ and $\Gamma'$ of
the same combinatorial graph, we say they are
\emph{$(\delta,\theta)$--close} if there exists a homeomorphism
between them such that corresponding points are within
distance~$\delta$ of each other, and corresponding
tangent vectors are within angle~$\theta$ of each other almost
everywhere.
\end{defn}

Remember~\cite{Sul-FTC} that any curve or graph~$\Gamma$ of
finite total curvature has well-defined one-sided tangent vectors everywhere;
these are equal and opposite except at countably many \emph{corners}
of~$\Gamma$ (including of course the vertices of the graph).
In the definition of $(\delta,\theta)$--close
it would be equivalent---in the case of \FTC/ or
piecewise~$C^1$ graphs---to require that the one-sided
tangent vectors be \emph{everywhere} within angle~$\theta$.

\section{Isotopy for thick knots}

Our first {isotopy} result looks at thick curves (that is,
embedded $C^{1,1}$ curves).
Here we can get close to optimal bounds on the~$\delta$
and~$\theta$ needed to conclude that nearby curves are isotopic.

The \emph{thickness}~$\tau(K)$\index{thickness (of curve)}
of a space curve~$K$ is
defined~\cite{gm} to be twice the infimal radius of
circles through any three distinct points of~$K$.
A link is~$C^{1,1}$ (that is, $C^1$ with Lipschitz tangent vector,
or equivalently with a weak curvature bound)
if and only if it has positive thickness~\cite{dsCKS2}.
Of course when~$K$ is~$C^1$, we can define normal tubes around~$K$,
and then $\tau(K)$ is the supremal diameter of such a tube
that remains embedded.  (We note that in the existing literature
thickness is sometimes defined to be the radius rather than diameter
of this thick tube.)  Points inside the tube have a unique
nearest neighbor on~$K$: the tube's radius is also the \emph{reach}
of~$K$ in the sense of Federer~\cite{federer}.

Fix a link~$K$ of thickness $\tau>0$.  Any other link which
follows~$K$ once around within its thick tube, transverse
to the normal disks, will be isotopic to~$K$.  This might seem
to say that any link $(\tau/2,\pi/2)$--close to~$K$ is isotopic,
but this would only work if we knew the correspondence used
in the definition of closeness were the closest-point projection to~$K$.
We do see that $(\delta,\theta)$--closeness is no guarantee
of isotopy if $\theta>\pi/2$ (allowing local knotting) or
if $\delta>\tau/2$ (allowing global strand passage).
So the following proposition (which extends Lemma~6.1 from~\cite{dsDDS},
with essentially the same proof) is certainly close to optimal.

(We note, however, that $\delta>\tau/2$ only allows global strand
passage in the case where the thickness of~$K$ is controlled
by self-distance rather than curvature.  This caveat explains why
we can later get some analogs of the following proposition even
for curves whose thickness is zero because they have corners.)

\begin{proposition}\label{pr:c11isot}
Fix a $C^{1,1}$ link~$K$ of thickness $\tau>0$.  For any
$\delta<\tau/4$, set $\theta := \pi/2-2\arcsin(2\delta/\tau)$.
Then any (rectifiable) link~$K'$ which is $(\delta,\theta)$--close
to~$K$ is ambient isotopic; the isotopy can be chosen to
move no point a distance more than~$\delta$.
\end{proposition}

\begin{proof}
For simplicity, rescale so that $\tau=1$, meaning the curvature
of~$K$ is bounded by~$2$.
Suppose $p'$ and~$p$ are corresponding points on~$K$ and~$K'$,
so $|p'-p|<\delta$.  Let $p_0$ be the (unique) closest point
on~$K$ to~$p'$.  Then of course $|p'-p_0|\le|p'-p|<\delta$,
so $|p-p_0|<2\delta$.  Standard results
on the geometry of thick curves~\cite[Lem.~3.1]{dsDDS}
show that the arclength of~$K$ between~$p$ and~$p_0$
is at most $\arcsin(2\delta)$, so by the curvature bound, the angle
between the tangents $T_p$ and~$T_{p_0}$ at these two points
is at most $2\arcsin(2\delta)$.
The definition of $(\delta,\theta)$--close means that
the tangent vector (if any) to~$K'$ at~$p'$ makes
angle less than~$\theta$ with~$T_p$,
so by definition of~$\theta$ it makes
angle less than~$\pi/2$ with~$T_{p_0}$.

This shows that $K'$ is transverse to the
foliation of the thick tube around~$K$ by normal disks.
We now construct the isotopy from~$K'$ to~$K$ as the union of
isotopies in these disks; on each disk we move $p'$ to~$p_0$,
coning this outwards to the fixed boundary.  No point in the
disk moves further than $p'$ does, and this is less than~$\delta$.
\end{proof}

\section{Isotopy for graphs of finite total curvature}

When generalizing our isotopy results to knotted graphs,
we must allow corners at the vertices of the graph, so it is
natural to consider \FTC/ graphs with possibly additional
corners along the arcs.  Our main theorem will again
show that sufficiently close curves are isotopic.

Any curve with a corner has zero thickness~$\tau=0$,
so it is perhaps surprising that an analog of \prop{c11isot}
still holds.  Perhaps the maximum possible
value of~$\delta$ here could be taken as a different
notion thickness of the curve; we do not explore this
idea, but content ourselves with merely proving the
existence of some positive~$\delta$.

A preliminary combing lemma allows the main theorem to
treat vertices on knotted graphs as well as sharp
corners along a single strand.  It was inspired
by Alexander and Bishop's correction~\cite{dsAlxBsh:FM}
to Milnor's treatment of corners of angle~$\pi$.
This case---where two or more strands leave a corner
in the same direction, allowing infinite winding or
braiding---is the most difficult, and seems to be the
obstacle to creating an isotopy which is $(\delta,\theta)$--small in
the sense described in Section~\ref{sec:small}.
But this case does not need to be treated separately in our proof.

This combing lemma is an example of the ``\ix{Alexander trick}''
(compare \cite[Prob.~3.2.10]{Thurston-3dGT}).
In general the trick shows that any homeomorphism of the unit $d$-ball
which fixes the boundary is isotopic to the identity.

\begin{figure} \centering
\includegraphics[width=.25\textwidth]{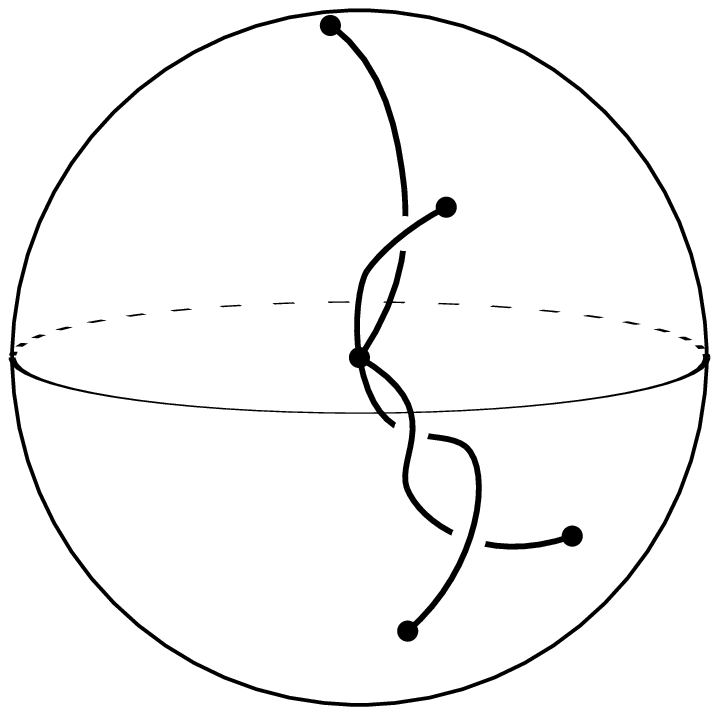} \qquad
\includegraphics[width=.25\textwidth]{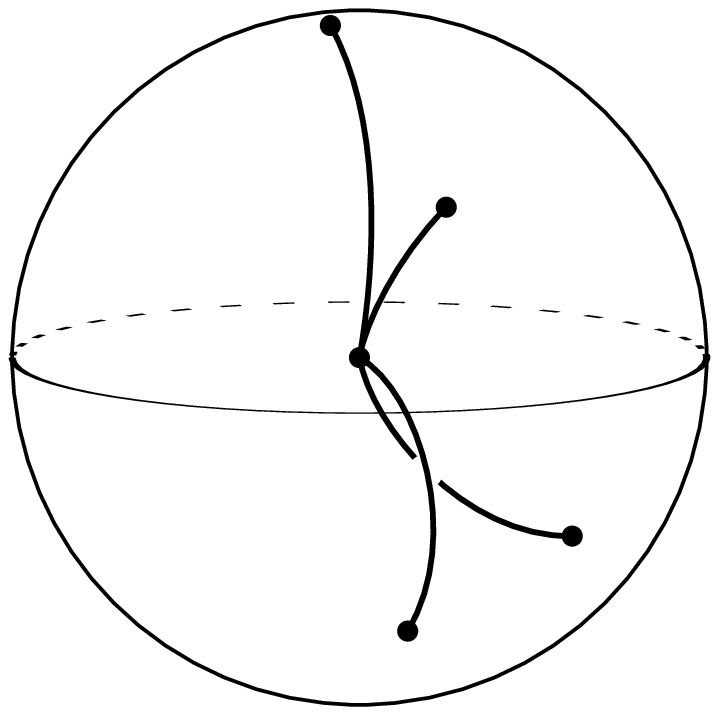} \qquad
\includegraphics[width=.25\textwidth]{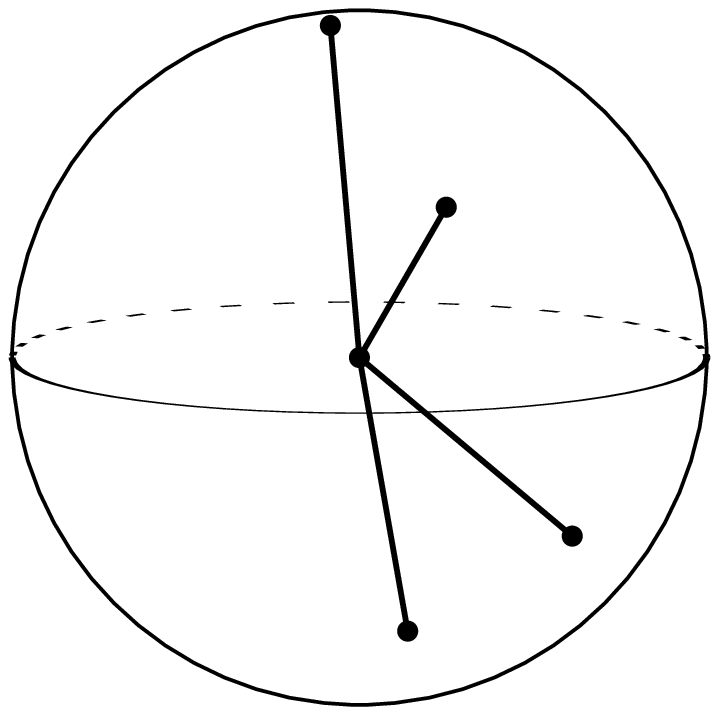}
\caption[Combing isotopy]{The isotopy of \lem{strball}
straightens strands in a ball by combing them inwards.
Only the outer half of the time~$0$ pattern (left)
is still seen at time~$\half$ (center);
at time~$1$ (right) we have straight radii.}\label{fig:comb}
\end{figure}

\begin{lemma}\label{lem:strball}
Let $B$ be a round ball centered at~$p$,
and suppose graphs $\Gamma$ and~$\Gamma'$
each consist of $n$ labeled arcs starting at~$p$ and proceeding
out to~$\bdry B$ transverse to the nested spheres around~$p$.
Then $\Gamma$ and~$\Gamma'$
are ambient isotopic.  In fact, any isotopy of~$\bdry B$ which takes the
$n$ points of~$\Gamma$ to those of~$\Gamma'$ can be extended to an isotopy
of the whole ball.
\end{lemma}
\begin{proof}
It suffices to show that $\Gamma$ is isotopic to a configuration
with $n$ straight radial segments to the same $n$ points on $\bdry B$.
For then the same is true of $\Gamma'$, and any two straight
configurations are clearly isotopic
(for instance by extending radially any given boundary isotopy).

To straighten $\Gamma$, we comb the $n$ strands inwards from~$\bdry B$,
which we identify with the unit sphere.
For $0<\lambda\le1$, let $p_i(\lambda)$ be the $n$ points of
intersection of~$\Gamma$ with the sphere $\lambda\bdry B$.
Start with~$f_1$ being the identity map, and define a family of maps
$f_\lambda$ from the unit sphere to itself, continuous in~$\lambda$,
such that $f_\lambda$ takes $p_i(1)$ to $p_i(\lambda)/\lambda$ for each~$i$.
The combing isotopy, at each time~$t$, maps each concentric sphere
$\lambda\bdry B$ to itself, using the map
$f_{\lambda+(1-\lambda)t}\after f_\lambda^{-1}$.
Thus at time~$t$ on the sphere of size~$\lambda$, we see the same
picture that was initially on the sphere of size $\lambda+(1-\lambda)t$;
see \figr{comb}.
It follows that at time $t=1$ the strands are all straight.
\end{proof}

As we have mentioned, if two (or more) strands
of~$\Gamma$ leave~$p$ with a common tangent, then they might
twist (or braid) around each other infinitely often.
But outside any sphere $\lambda\bdry B$
they can twist only finitely many times, explaining why
the isotopy constructed in the lemma has no trouble combing this out.

\begin{theorem}\label{thm:ftcgraph}
Suppose $\Gamma$ is a knotted graph of finite total curvature and $\eps>0$
is given.  Then there exists $\delta>0$ such that any
(rectifiable) graph~$\Gamma'$ which is $(\delta,\cutoff)$--close
to~$\Gamma$ is ambient isotopic to~$\Gamma$, via an isotopy
which moves no point by more than~$\eps$.
\end{theorem}
\begin{proof}
We begin by selecting a finite number of points~$p_j$ on~$\Gamma$
(including all its vertices) such that these points divide~$\Gamma$
into arcs~$\alpha_k$ each of total curvature less than~$\pi/8$.
(Note that any corner in~$\Gamma$ of turning angle at least~$\pi/8$
must be included among the~$p_j$.)

Let $r_1$ be the minimum
distance between any two arcs~$\alpha_k$ which are not incident to a
common~$p_j$ (or the minimum distance between points~$p_j$, if this
is smaller). Set $r_2 := \min(r_1/2,\eps/2)$.
Consider disjoint open balls~$B_j$ of radius~$r_2$ centered at the~$p_j$.
Each arc~$\alpha_k$ leaving~$p_j$
proceeds monotonically outwards to the boundary of~$B_j$
(since its curvature is too small to double back).
Also, $B_j$ contains no other arcs (since $r_2$ was chosen small enough);
indeed no other arcs come within distance~$r_2$ of~$B_j$.

Note that $\Gamma\setm\bigcup B_j$ is a compact union of disjoint arcs 
$\beta_k\subset\alpha_k$.  Let $r_3$ be the minimum
distance between any two of these arcs~$\beta_k$.
Considering endpoints of various~$\beta_k$ on the boundary of a single~$B_j$,
we note that $r_3\le 2r_2<\eps$.

For each~$\beta_k$, we construct a tube~$T_k$ around it as follows.
Suppose~$p_j$ and~$p_{j'}$ are the two endpoints of $\alpha_k\supset\beta_k$.
Foliate $\R^3\setm (B_j\cup B_{j'})$ by spheres of radius at least $r_2$
(and one plane, bisecting $\segment{p_jp_{j'}}$) in the obvious smooth way.
(Their signed curvatures vary linearly along $\segment{p_jp_{j'}}$.)
Note that because the arc $\alpha_k$ has total curvature less than~$\pi/8$,
by \cite[Lem.~2.3]{Sul-FTC} it is contained in a spindle of revolution
from~$p$ to~$p'$, of curvature~$\pi/4$.  Within this spindle,
the normal vectors to the foliating spheres stay within angle~$\pi/8$
of $\segment{p_jp_{j'}}$, as do the (one-sided) tangent vectors to~$\alpha_k$.
Thus $\beta_k$ is within angle $\pi/4$ of being normal to the foliation;
in particular it is transverse.  Finally, we set $r_4:=r_3/6$ and, for
each point $q\in\beta_k$, we consider the foliating sphere through~$q$
and in particular its intersection with~$B_{r_4}(q)$,
a (slightly curved) disk~$D_q$.
The tubular neighborhood~$T_k$ is the union of these disks
over all $q\in\beta_k$.  Since $r_4<r_2/3$, the normal vector
along each $D_q$ stays within angle $2\arcsin\tfrac16<\pi/8$
of that at its center~$q$, hence within $\pi/4$ of $\segment{p_jp_{j'}}$.

\begin{figure} \centering
\includegraphics[width=.35\textwidth]{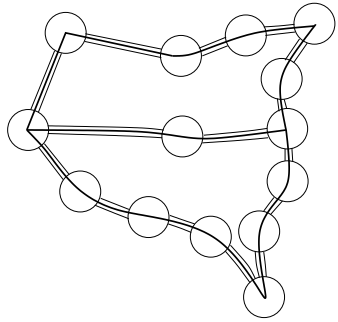}
\caption[Balls and tubes]{An \FTC/ graph has a neighborhood
consisting of balls and tubes, inside of which we perform
our isotopy to any nearby graph.}\label{fig:theta}
\end{figure}

We have thus found a neighborhood~$N$ of~$\Gamma$ which, as shown
in \figr{theta}, is foliated
by round spheres near the~$p_j$ and almost flat disks along the~$\beta_k$;
the graph~$\Gamma$ is transverse to this foliation~$F$.

Now take $\delta := r_4/3<r_2/9$ and consider some~$\Gamma'$ which is
$(\delta,\cutoff)$--close to~$\Gamma$.
We will see that~$\Gamma'$ lies within our neighborhood~$N$
and is also transverse to the foliation~$F$, except perhaps near the~$p_j$.
To find an ambient isotopy from~$\Gamma'$ to~$\Gamma$, we proceed in
two stages.

First we find the points $p'_j\in\Gamma'$ corresponding to $p_j\in\Gamma$,
and construct an isotopy $I$ supported in the $B_j$ as follows:
the map moves each point along a straight line;
$p'_j$ moves to $p_j$, each boundary point $q\in\d B_j$ is fixed,
and the map is linear on each segment $\segment{p_jq}$.
Because $\Gamma'$ was $(\delta,\cutoff)$--close to~$\Gamma$,
and because $p'_j$ was so close to $p_j$ compared to the radius~$r_2$ of~$B_j$,
the resulting $\Gamma'':=I_1(\Gamma')$ at time~$1$ has radial strands
within each~$B_j$, so we will be able to use the combing lemma.

In fact, we claim that~$\Gamma''$ is contained in~$N$ and
transverse to the foliation~$F$.  The original correspondence
between~$\Gamma'$ and~$\Gamma$ (showing they were
$(\delta,\cutoff)$--close) did not pair points in
the same leaf of the foliation.  But outside the~$B_j$,
for each point~$q'$ on the arc of $\Gamma''=\Gamma'$ corresponding
to~$\beta_k$, consider the leaf of the foliation it lies in
and the unique point $q_0\in\beta_k$ on the same leaf.
The distance $\big|q'-q_0\big|$ might be slightly bigger
than the distance $|q'-q|<\delta$ in the original pairing.
However, since the angles of the curves
and the foliation are so well controlled, it is certainly
well less than~$2\delta<r_4$; the tangent directions at~$q$
and~$q_0$ are certainly within $\pi/8$ of each other.
That is, even under the new correspondance ($q'$~to~$q_0$)
these arcs of~$\Gamma'$ and~$\Gamma$ are $(2\delta,\pi/4)$--close.
This proves the claim.

Second, we construct an isotopy~$J$ which is supported in~$N$ and
which preserves each leaf of the foliation~$F$.  On the disks in the
tubes~$T_k$ this is easy: we have to move one given point in the disk
to the center, and can choose to do this in a continuous way.
On the boundary $\d B_j$ we define $J$ to match these motions in
the disks and fix the rest of the sphere.  Finally, we use \lem{strball}
to fill in the isotopy~$J$ in the interior of the balls.

By construction, the overall isotopy~$I$ followed by~$J$ moves
points less than~$\eps$, as desired: within each ball~$B_j$ we
have little control over the details but certainly each point
moves less than the diameter $2r_2<\eps$; within~$T_k$ each
point moves less than $r_4\ll\eps$.
\end{proof}

Even though we have emphasized the class of \FTC/ curves,
we note that a very similar proof could be given for
the case when $\Gamma$ is piecewise~$C^1$.  We would
simply choose the points~$p_j$ such that the
intervening arcs~$\alpha_k$, while potentially of infinite
total curvature, had tangent vectors staying within
angle~$\pi/8$ of the vector between their endpoints.
Compare the beginning of the proof of \prop{smalliso} below.

As we noted in the introduction, we can recover the
following result: given any \FTC/ (or~$C^1$) link,
any sufficiently fine \ix{inscribed polygon} is isotopic.
When the link is~$C^1$, this is an immediate corollary,
sine the polygon will be $(\delta,\theta)$--close.
For an \FTC/ link~$K$ with corners, even very finely
inscribed polygons will typically cut those corners
and thus deviate by more than angle~$\cutoff$ from
the tangent directions of~$K$.  But, using \lem{strball}
near those sharp corners, we see immediately that the polygon
is isotopic to one that does use the corner as a vertex,
and thus to~$K$.

\section{Tame and locally flat links and graphs}\label{sec:tame}

Remember that a \ix{tame link} is one isotopic to a polygon.
Tame links clearly satisfy the following local condition.

\begin{defn}
A link $K$ is \dfn{locally flat} if each point $p\in K$ has
a neighborhood~$U$ such that $(U,U\cap K)$ is
an unknotted ball-arc pair.  (Unknotted means that the pair
is homeomorphic to a round ball with its diameter.)
\end{defn}

We can generalize the definitions to knotted graphs.

\begin{defn}
A knotted graph $\Gamma$ is \emph{tame} if it is isotopic
to a polygonal embedding; it is \emph{locally flat} if
each point $p\in \Gamma$ has a neighborhood~$U$ such that
$(U,U\cap \Gamma)$ is homeomorphic to a standard model.
At a $k$-fold vertex this standard model
is a round ball with $k$ of its radii;
along an arc of $\Gamma$ it is an unknotted ball-arc pair.
\end{defn}

Again, tame graphs are clearly locally flat.
In the mid-1950s, Bing~\cite{Bing} and Moise~\cite{Moise} independently
proved the much more difficult converse: any locally flat graph is tame.

Our main theorem showed that an \FTC/ graph is isotopic to
any other nearby graph.  This implies that \FTC/ graphs are
tame: as long as we include sharp corners among the polygon
vertices then sufficiently fine inscribed polygons will be
$(\delta,\cutoff)$--close.

Here we note that it is much easier to prove directly that
\FTC/ graphs are locally flat.

\begin{proposition}\label{pr:locflat}
Graphs of finite total curvature are locally flat.
\end{proposition}
\begin{proof}
Let $p_0\in\Gamma$ be any point along a graph of finite total curvature.
Repeat the construction at the beginning of the proof of \thm{ftcgraph},
but include~$p_0$ among the~$p_j$.  This gives a ball~$B_0$
around~$p_0$ containing only $k$ radial arcs of~$\Gamma$.
By \lem{strball} this ball is homeomorphic to any other such
ball, in particular to the standard model.
\end{proof}


\section{Applications to essential arcs}

One strand of work in \ix{geometric knot theory} attempts
to show that knotted curves are geometrically more
complex than unknots.  One of the earliest results in this direction is the
F{\'a}ry/Milnor theorem, which says a knotted
curve has more than twice the total curvature
of a round circle.  Interesting measures of
geometric complexity include ropelength~\cite{dsCKS2}
and distortion~\cite{dsKS-disto}.

\begin{figure}
\begin{overpic}[scale=.4]{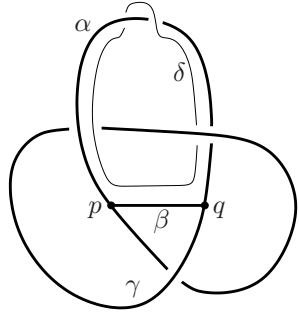}
\put(21.5,90){$\alpha$}
\put(26,31.5){$p$}
\put(66,31.5){$q$}
\put(47,26.5){$\beta$}
\put(38,6){$\gamma$}
\put(53.2,74){$\delta$}
\end{overpic}
\caption[An essential arc]
{In this knotted $\Theta$-graph the
ordered triple $(\alpha,\beta,\gamma)$ is essential.
The curve~$\delta$ is the parallel to $\alpha\cup\beta$
having linking number zero with $\alpha\cup\gamma$.
The fact that it is homotopically nontrivial in the knot
complement $\R^3\setm(\alpha\cup\gamma)$ is the obstruction
to the existence of a spanning disk.
In this illustration, $\beta$ is the straight segment~$\segment{pq}$, so
we equally say that the arc~$\alpha$ of the knot $\alpha\cup\gamma$
is essential.}
\label{fig:essdef} \end{figure}

We have recently obtained new lower bounds
for both ropelength~\cite{dsDDS} and distortion~\cite{dsDS-disto}
of knotted curves using the notion of essential arcs.
Generically, a knot~$K$ together with one of
its chords~$\segment{pq}$ forms a $\Theta$-graph in space;
being essential is a topological feature of this knotted graph.
The following definition, introduced in~\cite{dsDenne,dsDDS} and illustrated
in \figr{essdef}, is an extension of ideas of Kuperberg~\cite{Kup}.

\begin{defn}
Suppose $\alpha$, $\beta$ and~$\gamma$ are three disjoint simple arcs
from~$p$ to~$q$, forming a knotted $\Theta$-graph in~$\R^3$.  We say that
$(\alpha,\beta,\gamma)$ is \emph{inessential} if there is a disk~$D$
bounded by the knot $\alpha\cup\beta$ and having no interior intersections
with the knot $\alpha\cup\gamma$.  (We allow self-intersections of~$D$,
and interior intersections with~$\beta$; the latter are necessary if
$\alpha\cup\beta$ is knotted.)

Now suppose $p$ and~$q$ are two points along a knot~$K$,
dividing it into complementary subarcs $\arc{pq}$ and~$\arc{qp}$.
We say $\arc{pq}$ is \emph{essential} in~$K$ if for every $\eps>0$ there
exists some $\eps$-perturbation~$S$ of~$\segment{pq}$ (with endpoints fixed)
such that $K\cup S$ is an embedded $\Theta$-graph in which
$(\arc{pq},S,\arc{qp})$ is essential.
\end{defn}

\index{essential arc (of knot)}%
Allowing the $\eps$-perturbation ensures that the set of
essential arcs is closed in the set of all subarcs of~$K$;  it also
lets us handle the case when~$S$ intersects~$K$.  We could allow the
perturbation only in case of such intersections; applying an
argument like \lem{strball} in the case of an embedded~$\Theta$
(no intersections) shows this definition would be equivalent,
at least when the knot~$K$ is \FTC/.
We require only that the perturbation be small in a pointwise sense.
Thus $S$ could be locally knotted, but since the disk in the definition
can intersect $\beta=S$ this is irrelevant.

This definition of essential in terms of $\Theta$-graphs
was what led us to consider knotted graphs and to look
for a result like \thm{ftcgraph}.  We now apply that theorem
to show that essential arcs are preserved in limits.
(In \cite[Prop.~6.2]{dsDDS} we gave a similar result, but
since we did not yet have \thm{ftcgraph}, we had to
restrict to $C^{1,1}$ knots only.  The proof below does
closely follow the proof there, \thm{ftcgraph} being
the main new ingredient.) 

\begin{proposition}
Suppose a sequence of rectifiable knots~$K_i$ converge to an \FTC/ limit~$K$
in the sense that $K_i$ is $(\delta_i,\pi/10)$--close to~$K$
with $\delta_i\to0$.
Suppose that the~$K_i$ have essential subarcs $\arc{p_iq_i}$ with
$p_i\rightarrow p$ and $q_i\rightarrow q$.
Then~$\arc{pq}$ is an essential subarc of~$K$.
\end{proposition}

\begin{proof}
We can reduce to the case $p_i=p$, $q_i=q$ (but with $\pi/10$
above replaced by $\cutoff$) by appyling
euclidean similarities (approaching the identity) to the $K_i$.

Given any $\epsilon>0$, we prove there is an $2\epsilon$-perturbation~$S$
of~$\segment{pq}$ for which $(\arc{pq},S,\arc{qp})$ is essential.
Then by definition $\arc{pq}$ is essential.

We want to find an ambient isotopy~$I$ from some~$K_i$ to~$K$ which moves
points less than~$\eps$.  So apply \thm{ftcgraph} to~$K$ to determine
$\delta>0$ and then pick~$i$ large enough so that
the knot~$K_i$ is $(\delta,\cutoff)$--close to~$K$.
The ambient isotopy~$I$ guaranteed by \thm{ftcgraph}, with $K=I(K_i)$,
moves points less than~$\eps$, as desired.

Since $\arc{p_iq_i}\subset K_i$ is essential,
by definition, we can find an $\epsilon$-perturbation $S_i$
of~$\segment{p_iq_i}$ such that
$(\arc{p_iq_i},S_i,\arc{q_i,p_i})$ is essential.
Setting $S:=I(S_i)$, this is the desired $2\epsilon$-perturbation
of $\segment{pq}$.  By definition, the $\Theta$-graph $K_i\cup S_i$
is isotopic via~$I$ to $K\cup S$,
so in the latter $(\arc{pq},S,\arc{qp})$ is also essential.
\end{proof}

\section{Small isotopies in a stronger sense}\label{sec:small}

The isotopy constructed in \prop{c11isot} is small in a pointwise sense.
But we have required our curves to be close in a stronger sense.
We will say that a $C^1$~isotopy~$I$, or more precisely its time-$1$
map~$I_1$, is \emph{$(\delta,\theta)$--small} if $I_1$ moves no
point distance more than~$\delta$ and its derivative turns no
tangent vector more than angle~$\theta$.
Note that under such a small isotopy, every rectifiable graph~$\Gamma$
is $(\delta,\theta)$--close to its image~$I_1(\Gamma)$.
(Although we could extend the definition of ``small''
for certain nonsmooth isotopies, we do not pursue this idea here.
It remains unclear if the ball isotopy of \lem{strball} could be made
small in some such sense.)

\begin{lemma}\label{lem:lip-isot}
Suppose an ambient isotopy is given for $t\in[0,1]$
by $I_t(p) = p+tf(p)$, where $f:\R^3\to\R^3$ is smooth,
bounded in norm by~$\delta$, and $\lambda$-Lipschitz.
Then $I$ is $(\delta,\arctan\lambda)$--small.
\end{lemma}
\begin{proof}
Clearly each point $p\in\Gamma$ is within distance~$\delta$ of
$p':=I_1(p)$.  If~$v$ is the unit tangent vector to~$\Gamma$
at~$p$, then its image, the tangent vector to~$\Gamma'$
at~$p'$, is $v+\d_v f$, where $|\d_v f|<\lambda$ by the Lipschitz
condition.  This differs in angle from~$v$ by at most $\arctan\lambda$.
\end{proof}

The isotopy we built disk-by-disk in \prop{c11isot}
is not necessarily smooth, and does not necessarily
satisfy the hypotheses of this lemma.
In particular, since~$K$ might not be piecewise~$C^2$, its normal
disks might not form a piecewise~$C^1$ foliation, and the resulting
isotopy would not even have one-sided derivatives everywhere.

So we will use a different construction, not using the thick tube.
Here, as in \thm{ftcgraph}, we abandon any attempt to get optimal
constants, but the starting curve can be any $C^1$ curve,
as in Crowell and Fox's construction \cite[App.~I]{CFox}
of an isotopic inscribed polygon.
We begin with a lemma which captures what we feel is the
essential ingredient of that argument.

\begin{lemma}\label{lem:c1-unif}
Given a $C^1$~link~$K$ and an angle $\theta>0$,
there is some $\ell>0$ such that along any subarc~$\alpha$
of~$K$ of length at most~$\ell$, the tangent vector stays
within angle~$\theta$ of the chord vector~$v$ connecting
the endpoints of~$\alpha$.
\end{lemma}
\begin{proof}
Because~$K$ is compact, its continuous tangent vector is
in fact uniformly continuous.  Noting that~$v$ is an
average of the tangent vectors along~$\alpha$, the
result follows immediately.
\end{proof}

This implies that the distortion
(the supremal arc/chord ratio) of any~$C^1$ curve is bounded:
by the lemma the ratio approaches~$1$ for short arcs,
so the supremum is achieved.

\begin{proposition}\label{pr:smalliso}
Given a $C^1$~link~$K$, we can find $\eps>0$ such that
the following holds.  For any $\delta<\eps$ and $\theta<\pi/6$,
and any $C^1$~link~$K'$ which is $(\delta,\theta)$--close to~$K$,
there is a $(2\delta,2\theta)$--small ambient isotopy
taking~$K'$ to~$K$.
\end{proposition}
\begin{proof}[Sketch of proof]
Fix $\phi<5^\circ$, and find $\ell>0$ as in \lem{c1-unif}.
Now set $\tau := \min |x-y|$, where the minimum is taken
over all points $x,y\in K$ not connected by a subarc of
length less than~$\ell$.  (Note that $\tau<\ell$.)

Along each component of~$K$, pick $r\in(\tau/50,\tau/40)$ and
place points~$p_j$ spaced equally at arclength~$r$.
The resulting inscribed polygon~$P$ has edgelengths
in the interval $[r\cos\phi,r]$
and is $(r\sin\phi,\phi)$--close to~$K$,
with turning angles at the~$p_j$ less than~$2\phi$.
Rounding off the corners in a suitable fashion,
we obtain a $C^2$~curve~$L$ whose radius of curvature is
never less than $5r$, and which is
$(d,\phi)$--close to~$P$ for $d\le(\sec\phi-1)r/2<r/500$.
We see that $L$ is $(\tau/400,2\phi)$--close to~$K$,
and has radius of curvature at least~$\tau/10$.

We claim that the normal tube around~$L$ of diameter~$\tau/5$
is embedded.  This is true locally by the curvature bound.
By standard results on thickness (cf.~\cite{dsCKS2,dsDDS})
it can then only fail if there is a doubly-critical pair
on~$L$ at distance less than~$\tau/5$: a pair whose chord
is perpendicular to~$L$ at both ends.  Then there would
be a nearby pair on~$K$, still at distance less than~$\tau/4$.
Since an arc of~$L$ connecting the doubly critical pair
must turn at least~$\pi$, the pair on $X$ cannot be connected
by an arc of length less than~$\ell$.
Thus we have a contradiction to the choice of~$\tau$.

Now choose $\eps:=\tau/100$.
For any $\delta<\eps$ and $\theta<\pi/6$,
if $K'$ is $(\delta,\theta)$--close to~$K$, then
it is easy to check that both are $C^1$, transverse
sections in the tube of radius~$\tau/50$ around~$L$.
We construct an isotopy between them, supported in
the embedded tube of radius~$\tau/10$ around~$L$.
On each normal disk, we move a subdisk of
radius~$\tau/25$ rigidly, to move its center from~$K$ to~$K'$.
In the outer part of the thick tube, we extend in
such a way that the isotopy is $C^1$ throughout space.
The isotopy moves no points more than the distance
between~$K$ and~$K'$ in the normal disk.  This may be
somewhat more than~$\delta$, since the original pairing
between the two links may have been different, but is
certainly less than~$2\delta$.

To check that this isotopy turns tangent vectors at most~$2\theta$,
it is easiest to apply \lem{lip-isot}.  The Lipschitz
constant~$\lambda$ is basically determined by needing
to turn~$K'$ by approximately angle~$\theta$ to match~$K$.
Because the normal disks get slightly closer to one another
as we move off~$L$ in the direction of its curvature,
and since we have smoothed outside the thinner tube,
the constant goes up slightly, but not enough to
violate the~$2\theta$. 
\end{proof}
\index{isotopy|)}


\subsection*{Acknowledgments}
We extend our thanks to Stephanie Alexander, Dick Bishop
and Dylan Thurston for their interest and helpful conversations.



\providecommand{\bysame}{\leavevmode\hbox to3em{\hrulefill}\thinspace}

\end{document}